\documentclass[10pt]{article}

\usepackage[all]{xy}
\usepackage{textcomp}
\usepackage{amsmath,amssymb,eucal,eufrak}
\usepackage{graphicx}
\newtheorem{theorem}{Theorem}[section]
\newtheorem{lemma}[theorem]{Lemma}
\newenvironment{proof}{\par{\bf Proof\ }}{}

\newtheorem{definition}[theorem]{Definition}
\newtheorem{example}[theorem]{Example}

\newtheorem{proposition}[theorem]{Proposition}
\newtheorem{corollary}[theorem]{Corollary}
\newtheorem{remark}[theorem]{Remark}

\numberwithin{equation}{section}

\def\g{{\mathrm{g}}}  \def\h{{\mathrm{h}}}

 \def\rK{{\mathrm{K}}}

  \def\CC{{\mathbb{C}}}
  
  \def\II{{\mathbb{I}}}
 \def\KK{{\mathbb{K}}} 
 \def\NN{{\mathbb{N}}} 
  \def\RR{{\mathbb{R}}}
  \def\UU{{\mathbb{U}}}

 \def\bK{{\mathbf{K}}} \def\bL{{\mathbf{L}}}

 \def\cB{{\mathcal{B}}} 
\def\cD{{\mathcal{D}}} \def\cE{{\mathcal{E}}} \def\cF{{\mathcal{F}}}
\def\cG{{\mathcal{G}}} \def\cH{{\mathcal{H}}} 
 \def\cK{{\mathcal{K}}} \def\cL{{\mathcal{L}}}
\def\cM{{\mathcal{M}}}  
\def\cP{{\mathcal{P}}}  \def\cR{{\mathcal{R}}}
 \def\cT{{\mathcal{T}}} \def\cU{{\mathcal{U}}}
 \def\cW{{\mathcal{W}}} 
 
  \def\gC{{\mathfrak{C}}}
  \def\gF{{\mathfrak{F}}}
 \def\gH{{\mathfrak{H}}} 
 \def\gK{{\mathfrak{K}}}

\def\gV{{\mathfrak{V}}}

\begin{document}

\title{Invariant {K}re\u{\i}n subspaces, regular irreducibility and integral representations}
%\title{Invariant {K}re\u{\i}n subspaces and their kernels, integral
%representations}

\author{Xavier Mary \footnote{email: xavier.mary@ensae.fr}\\
\textit{\small Ensae - CREST, 3, avenue Pierre Larousse 92245 Malakoff Cedex, France}}

\date{}
\maketitle
\begin{abstract}
We study unitary representations of groups in Kre\u{\i}n spaces, irreducibility criteria and integral decompositions. Our main tool is the theory of Kre\u{\i}n subspaces and their
(reproducing) kernels and a variant of Choquet's theorem.
\end{abstract}


%%%%%%%%%%%%%%%%%%%%%%%%%%%%%%%%%%%%%%%%%%%%%%%%%%%%%%%%%%%%
%
%%%%%%%%%%%%%%%%%%%%%%%%%%%%%%
%
%
%\section*{Introduction}
%
%
%%%%%%%%%%%%%%%%%%%%%%%%%%%%%%
%
%%%%%%%%%%%%%%%%%%%%%%%%%%%%%%%%%%%%%%%%%%%%%%%%%%%%%%%%%%%%

%%%%%%%%%%%%%%%%%%%%%%%%%%%%%%
\section*{Introduction}
%%%%%%%%%%%%%%%%%%%%%%%%%%%%%%

The use of reproducing kernels methods in harmonic analysis and representation theory is now classical \cite{Schwartz64Lie}, \cite{Faraut99}, \cite{Pestman85}, \cite{Thomas05},
\cite{VanDijk2002}, and is one of the main applications of the theory of reproducing kernel Hilbert spaces
or Hilbert subspaces \cite{Aronszajn50}, \cite{Schwartz64}.\\

The present paper was motivated by the intuition that such methods could also apply to representations in indefinite inner product spaces. Hence we consider together two less known
extensions of the previous theories: on one hand, the theory of Hermitian subspaces and Kre\u{\i}n subspaces and their kernels \cite{Schwartz64}, \cite{Alpay91}, \cite{Sorjonen73}, and
on the other hand operator algebras and group representations in indefinite inner product spaces \cite{Naimark65}, \cite{Loginov69}, \cite{Kissin96}, \cite{Ota88}, \cite{Ota89}. Note
that these representations appear mainly in mathematical physics, following the use of indefinite metric spaces in the works of Dirac \cite{Dirac42}, Pauli \cite{Pauli43} and
Heisenberg \cite{Heisenberg66}. They are now critical issues in quantum electromagnetism and the Gupta-Bleuler triplet \cite{Araki85}, representations of CCR algebras
\cite{Mnatsakanova98}, \cite{Mnatsakanova03} or the QFT formalism \cite{Strocchi78} and the study of De Sitter spaces, \cite{Gazeau06}. Note also that since the pioneering work of
Pontryagin \cite{Pontryagin44}, the theory of linear operators in Kre\u{\i}n spaces
has been developed into a major branch of modern operator theory \cite{Iokhvidov60}, \cite{Langer82}, \cite{Alpay01}.\\

The paper is organized as follows: in section 1 we review some general facts about Kre\u{\i}n spaces, and define Kre\u{\i}n subspaces and their kernels. Section 2 is devoted to
indefinite representations, and discuss the links between invariant Kre\u{\i}n subspaces and their kernels. It gives also criteria for irreducibility and a variant of Schur's lemma.
Sections 3 and 4 then discuss the existence of a direct integral decomposition into irreducibles subspaces.

\section{Kre\u{\i}n spaces, Kre\u{\i}n subspaces and kernels}

\subsection{Kre\u{\i}n spaces}
 A Kre\u{\i}n space is an indefinite inner product space
$(\cK, \left[ \cdot,\,\cdot \right])$ (\textit{i.e.} the form $\left[ \cdot,\,\cdot \right]$ is sesquilinear and Hermitian) such that there there exists an automorphism $J$ of $\cK$
which squares to the identity (called fundamental symmetry or signature operator), $\langle x,\,y\rangle \equiv \left[ Jx,\,y \right]$ defines a positive definite inner product and
$(\cK, \langle \cdot,\,\cdot \rangle)$ is a Hilbert space. Equivalently, the indefinite inner product space $(\cK, \left[ \cdot,\,\cdot \right])$ is a Kre\u{\i}n space if there exist
an admissible (with respect to the inner product) hilbertian topology on $\cK$ that makes it an Hilbert space.

%Kre\u{\i}n spaces are named in honor of the Ukrainian
%mathematician Mark Grigorievich Kre\u{\i}n.

 The following subsets are defined in terms of the ``square
norm'' induced by the indefinite inner product:
\begin{itemize}
\item
$\bK_{+} \equiv \{ x \in \cK : \left[ x,\,x \right] > 0 \}$ is called the ``positive cone"; \item $\bK_{-} \equiv \{ x \in \cK : \left[ x,\,x \right] < 0 \}$ is called the ``negative
cone";
\item $\bK_{0} \equiv \{ x \in \cK : \left[ x,\,x \right] =0  \}$ is
called the ``neutral cone".
\end{itemize}

 A subspace $\cL \subset \cK$ lying within $\bK_{0}$ is
called a ''neutral subspace''.  Similarly, a subspace lying within $\bK_{+}$ ($\bK_{-}$) is called ''positive'' (''negative''). A subspace in any of the above categories may be called
''semi-definite'', and any subspace that is not semi-definite is called ''indefinite''.

 Any decomposition of the indefinite inner product space
$\cK$ into a pair of subspaces $\cK = \cK_+ \oplus \cK_-$ such that $\cK_+ \subset \bK_{+}\cup\{0\}$ and $\cK_- \subset \bK_{-}\cup\{0\}$ is called a ''fundamental decomposition'' of
$\cK$. $\cK_+$ equipped with the restriction of the bilinear form $[\cdot,\cdot]$ is then a Hilbert space, and $\cK_-$ the antispace of a Hilbert space $|\cK_-|$. To this fundamental
decomposition is associated a fundamental symmetry $J$ such that the scalar product $\langle x,\,y\rangle \equiv \left[ Jx,\,y \right]$ coincide with the scalar product of
$\cH=\cK_+\oplus |\cK_-|$.

 The positive definite inner product
$\langle\cdot,\,\cdot\rangle$ depends on the chosen fundamental
decomposition, which is, in general, not unique.  But %(corollary 1p~5 in \cite{Dritschel96})
(see \cite{Dritschel96}) any two fundamental symmetries $J$ and $J^\prime$ compatible with the same indefinite inner product on $\cK$ result in Hilbert spaces $|\rK|$ and
$|\rK^\prime|$ whose decompositions $|\rK|_\pm$ and
$|\rK^\prime|_\pm$ have equal dimensions. Moreover %(corollary 2 p~5 in \cite{Dritschel96})
they induce equivalent square norms hence a unique topology. This topology is admissible, and it is actually the Mackey topology defined by the bilinear pairing. All topological
notions in a Kre\u{\i}n space, like continuity or closedness of sets are understood with respect to this Hilbert space topology.

Orthogonality is a key issue in indefinite inner product spaces. Let $\cL$ be a subspace of $\cK$. The subspace $\cL^{[\perp]} \equiv \{ x \in \cK : \left[ x,\,y \right] = 0$ for all
$y \in \cL \}$ is called the orthogonal companion of $\cL$. If $J$ is a fundamental symmetry it is related to the (Hilbert) orthogonal by $\cL^{[\perp]}=(J\cL)^{\perp}$. $\cL^{0}
\equiv \cL \cap \cL^{[\perp]}$ is called the isotropic part of $\cL$. If $\cL^{0} = \{0\}$, $\cL$ is called non-degenerate. It is called regular (or a Kre\u{\i}n subspace) if it is
closed and a Kre\u{\i}n space with respect to the restriction of the indefinite inner product. This is equivalent to $\cL\oplus \cL^{[\perp]}=\cK$ (\cite{Dritschel96}) and this
relation is sometimes taken as a definition of regular subspaces.

 If $\cK$ and $\cH$ are Kre\u{\i}n spaces, then continuity
of operators is defined with respect to the Hilbert norm induced by any fundamental decomposition.

Any continuous (weakly continuous) operator $A$ has an adjoint $A^{[*]}$
%to each others
(with respect to the indefinite inner product) verifying
\begin{equation}
\forall k\in \cK,\, \forall h\in \cH, \qquad [Ak,h]_{\cH}=[k,A^{[*]}h]_{\cK}
\end{equation}
This adjoint is sometimes called $J$-adjoint (to emphasize the role of the fundamental symmetry $J$) and it is related to the classical (Hilbert) adjoint by $A^{[*]}=JA^*J$.

\subsection{Kre\u{\i}n subspace and kernels}
%From now on we fix a pair $(\cE,\cF)$ of spaces together with a
%non-degenerate sesquilinear form $\left\dot,\dot\right)_{(\cF,\cE)}$
%such that $\cE$ endowed with the Mackey topolgy induced by the form
%is weakly complete. If we are given a weakly complete locally convex
%space $\cE$, then we can take for $\cF$ is topological antidual

We first make some definitions on remarks on kernels.

\begin{definition}
$\,$\\ \vspace{-1cm}
\begin{enumerate}
\item Let $\cE$ be a l.c.s., $\cE'$ its topological dual and $\overline{\cE'}$ the conjugate space of its topological dual. Then according to L. Schwartz \cite{Schwartz64},
we call kernel any weakly continuous linear application $\varkappa:\overline{\cE'}\longrightarrow\cE$.
\item More generally, if $(\cE,\cF)$ is a semi-duality (a pair $(\cE,\cF)$ of spaces together with a
non-degenerate sesquilinear form $\left(.,.\right)_{(\cF,\cE)}$), we call kernel any weakly continuous linear application $\varkappa:\cF\longrightarrow\cE$.
\end{enumerate}

\end{definition}

Since $\cF$ is identified with the space of continuous semilinear forms on $\cE$, the adjoint of $\varkappa$ ($\varkappa^*$) is also a kernel of $(\cE,\cF)$. The kernel $\varkappa$ is
Hermitian if $\varkappa^*=\varkappa$. It is positive if
\begin{equation}\label{eq kern pos}\forall \varphi\in \cF, \left(\varphi,
\varkappa(\varphi)\right)_{(\cF,\cE)}\geq 0
\end{equation}
A positive kernel is Hermitian. We note the set of positive kernels $\bL^{+}(\cF,\cE)$, the set of Hermitian kernels $\bL^{h}(\cF,\cE)$.

It is crucial to note that the set of positive kernels $\bL^{+}(\cF,\cE)$ is a salient convex cone (the positivity conditions then defines a partial order on the set of positive
kernels: $H\leq K\iff K-H\geq 0$.)

The set of Hermitian kernels will usually be to large for our applications and we will mainly study bounded kernels. A Hermitian kernel $K$ is bounded by a positive kernel $H$ if
$H-K\geq 0$ and $H+K\geq 0$. In this case $H$ is called a majorant of $K$ and we say that $(K,H)$ is a bounded pair (of kernels). We note the set of bounded hermitian kernels
$\bL^{b}(\cF,\cE)$.

We say that positive kernels $K$ and $L$ are independent if the following
statement holds:\\ If $0\leq H\leq K$ and $0\leq H\leq L$ then $H=0$.\\
We will sometimes use the following notation: if $K$ and $L$ are positive (kernels $K+L=K\oplus L$ means that the kernels are independent.

\begin{definition}
$(K, H)$ is a Kolmogorov Hermitian pair (or minimal pair) of kernels if $K$ is a Hermitian kernel, $H$ a positive kernel that bounds $K$ and $H-K, H+K$ are independent.
\end{definition}

Second, we define Hilbert and Kre\u{\i}n subspaces of a semi-duality $(\cE,\cF)$ (equivalently of a l.c.s. $\cE$).

\begin{definition}%[Hilbert subspace, Hermitian subspace, Kre\u{\i}n subspace]
$\,$\\ \vspace{-1cm}
\begin{enumerate}
\item A Hilbert subspace $\cH$ of $(\cE,\cF)$ is a Hilbert space
continuously (for the Mackey topologies) included in $\cE$.
\item A Hermitian subspace $\cK$ of $(\cE,\cF)$ is the difference of two (disjoint) Hilbert
subspaces of $(\cE,\cF)$.
\item A Kre\u{\i}n subspace $\cK$ of $(\cE,\cF)$ is a Kre\u{\i}n space
continuously (for the Mackey topologies) included in $\cE$.
\end{enumerate}
\end{definition}

It si straightforward to see that the last two notions coincide (for details on Hermitian subspaces an Kre\u{\i}n subspaces, we refer to \cite{MaryPHD}). We note $Hilb((\cE,\cF))$ the
set of Hilbert subspaces of $(\cE,\cF)$ and $Krein((\cE,\cF))$ the set of Kre\u{\i}n subspaces of $(\cE,\cF)$.

We can now state the main results of the theory of Hilbert subspaces and Kre\u{\i}n subspaces:

\begin{theorem}\label{isomorphism}
Suppose $\cE$ is quasi-complete (for its Mackey topology). Then there is a bijection between $Hilb((\cE,\cF))$ and $\bL^{+}(\cF,\cE)$. Moreover, this bijection is an isomorphism of
convex cones.
\end{theorem}

 Any reader particularly interested by the isomorphism of
convex cone structure can read \cite{Schwartz64} p~159-161, where the proof is detailed. The kernel $H$ of a Hilbert subspace $\cH$ is characterized by the following equality:
\begin{equation}\label{eqH}
\forall \varphi\in \cF,\; \forall h\in \cH,\quad \langle H(\varphi),h\rangle_{\cH}=\left(\varphi, h\right))_{\cF,\cE}
\end{equation}

\begin{theorem}\label{lack isomorphism}
There is a surjective morphism of convex cones between $Krein((\cE,\cF))$ and $\bL^{b}(\cF,\cE)$. This is \textbf{not} an isomorphism in general.
\end{theorem}
The kernel $K$ of a Kre\u{\i}n subspace $\cK$ is characterized by an equation similar to \ref{eqH}:
\begin{equation}\label{eqK}
\forall \varphi\in \cF,\; \forall k\in \cK,\quad \left[ K(\varphi),k\right]_{\cK}=\left(\varphi, k\right))_{\cF,\cE}
\end{equation}
The existence of so-called kernels of multiplicity is the major difference between the two theories. It follows that all the constructions that rely uniquely on the kernel may fail, as
we will see studying integral of Kre\u{\i}n subspaces. To circumvent this flaw we introduce the notion of Kre\u{\i}n-Hilbert pairs.

\subsection{Kre\u{\i}n-Hilbert pairs}

\begin{definition}%[Kre\u{\i}n-Hilbert pair]
A pair $(\cK,\cH)$ is called a Kre\u{\i}n-Hilbert pair $\cK$ is a Kre\u{\i}n
space continuously included in the Hilbert space $\cH$.\\
It is called a closed Kre\u{\i}n-Hilbert pair  if $\cK$ is a Kre\u{\i}n
space, closed subspace of the Hilbert space $\cH$.\\
It is a fundamental (or minimal) pair if there exists two Hilbert spaces $\cH_+$ and $\cH_-$,
$$\cK=\cH_+\ominus \cH_- \text {and } \cH=\cH_+\oplus \cH_-$$
\end{definition}
 Of course any fundamental pair is closed. Any pair
$(\cK,\cH)$ define a (bounded) hermitian kernel $\chi:\cH\to\cH$. It is interesting at this point to note that even in this case the kernel $\chi$ may be of multiplicity (see
\cite{Schwartz64}, \cite{Hara92}, \cite{Langer2003}).

\begin{example}
Let $(\cK,J)$ be a Kre\u{\i}n space with fundamental symmetry $J$. Then $\cH=\cK$ endowed with the scalar product $\langle k,h\rangle_{\cH} =[k, Jh]_{\cK}$ is a a Hilbert space and
$(\cK,\cH)$ is a fundamental Kre\u{\i}n-Hilbert pair.
\end{example}

We note $KH((\cF,\cE))$ the set of minimal Kre\u{\i}n-Hilbert pairs of subspaces $(\cK,\cH)$, equivalently of pairs $(\cK,\cH)$ of a Kre\u{\i}n subspace and a Hilbert subspace of
$\cE$, such that their kernels $(K,H)$ form a minimal pair.

%\begin{theorem} FAUX !!!
%Let $(\cK,J)$ be a Kre\u{\i}n space with fundamental symmetry $J$, $(\cK,
%\cH)$ the associated Kre\u{\i}n-Hilbert pair. Let $\cM$ be a regular
%subspace of $\cK$, $P$ the associated (orthogonal) projection and
%$\cN=cM^{[\perp]}$ its orthogonal complement. Finally let $Q$ be the
%projection on $\cM$ perpendicularly to $\cN$. Then $Q-P$ and $Q+P$
%are positive operators.
%\end{theorem}

%\begin{proof}
%
%\end{proof}

\subsection{Image by a weakly continuous application}

 We suppose now we are given a second pair of spaces in
duality $(\mathfrak{E},\mathfrak{F})$. It is actually possible (\cite{Schwartz64}, \cite{Mary2005} and proofs therein) to define the image of a Kre\u{\i}n space $\cK$ by a weakly
continuous linear application $u:\cE\to \mathfrak{E}$
 by using orthogonal relations in the duality $\cK$, but this image is not in general a Kre\u{\i}n space.\\
We recall this construction hereafter. $\forall A\subset \cE,\; u_{|A}$ denotes the restriction of $u$ to the set $A$. We then define the following quotient space:
$$\cM=\left(\ker(u_{|\cK})^{[\perp]}/\ker(u_{|\cK})\right)$$
\begin{lemma}
The linear application $u_{|\cM}$ is well defined and injective, and $\forall \,(\dot{m},\dot{n})\in \cM^2$, the bilinear form $B(u_{|\cM}(\dot{n}),u_{|\cM}(\dot{m}))=[n,m]_{(\cK)}$
defines a indefinite inner product on the space $u_{|\cM}(\cM)$.
\end{lemma}

\begin{proof}
We have the following factorisation
$$u:\ker(u_{|\cK})^{[\perp]}\longrightarrow (\ker(u_{|\cK})^{[\perp]}/\ker(u_{|\cK})\overset{u_{|\cM}}{\longrightarrow} \mathfrak{E}$$
and  $u_{|\cM}$ is one-to-one.
Moreover the bilinear form $B:u_{|\cM}(\cM)\times u_{|\cM}(\cM)\longrightarrow \KK$ is well defined since:\\
$\forall (m_{1},m_{2})\in \dot{m},\, \forall (n_{1},n_{2})\in \dot{n},\; [m_{1}-m_{2},n_{1}-n_{2}]_{\cK}=0$.
\end{proof}

% The space $u(\cM)$ is not a Kre\u{\i}n space in general, but is
%is weakly included

 However, it may happen that this image is Mackey-complete
(for instance if $u$ is one-to-one, or more generally if $\ker(u)$ is regular). In this case the space $u_{|\cM}(\cM)$ is actually a Kre\u{\i}n space continuously embedded in
$\mathfrak{E}$, and we can compute its kernel:
\begin{theorem}\label{image}
$\,$\\
 If $u_{|\cM}(\cM)$ is a Kre\u{\i}n space, then it is a Kre\u{\i}n
subspace of $\mathfrak{E}$. Its kernel is $u\circ\varkappa\circ ^{t}u$.
\end{theorem}

\subsection{Integral of Kre\u{\i}n subspaces: the neccessity of Kre\u{\i}n-Hilbert pairs}
The theory of direct integral of Hilbert spaces is well known, as is the theory of integral of Hilbert subspaces (\cite{Schwartz64}, \cite{Thomas05}), and poses no difficulties. This
is not the case for Kre\u{\i}n spaces, where it is actually not possible to define directly the direct integral of Kre\u{\i}n spaces as the following example shows:

\begin{example}\label{ex 1}
Define on $\RR^2$ the following inner products :
\begin{eqnarray*}
\left[\left(
       \begin{array}{c}
         x_1 \\
         y_1 \\
       \end{array}
     \right)
|\left(
       \begin{array}{c}
         x_2 \\
         y_2 \\
       \end{array}
     \right)\right]&=x_1y_2+x_2 y_1\\
     \langle\left(
       \begin{array}{c}
         x_1 \\
         y_1 \\
       \end{array}
     \right)|\left(
       \begin{array}{c}
         x_2 \\
         y_2 \\
       \end{array}
     \right)\rangle_n&=\frac{1}{n^2}x_1x_2+n^2 y_1 y_2
\end{eqnarray*}
Then $\left(\cK=\RR^2,\left[.|.\right]\right)$ is a Kre\u{\i}n space, and
$$J_n=\left(
       \begin{array}{cc}
         0 &  n^2 \\
         \frac{1}{n^2} & 0 \\
       \end{array}
     \right),\; n\in \NN^*$$

are fundamental symmetries associated with the spaces
$\left(\cH_n=\RR^2,\langle.|.\rangle_n\right)$.\\

We want to give a sense to $\gK=\displaystyle\bigoplus_{n=1}^{+\infty}\cK$.
 Fact is that there are many possible interpretations of this
 space in terms of Kre\u{\i}n space. For instance $$\gK_1=\gH_1=\bigoplus_{n=1}^{+\infty}\cH_1=\{h=k_1\oplus k_2\oplus  ...,\; \sum_{n=1}^{\infty} ||k_n||^2_1<+\infty\}$$ as a vector space, with inner product
 $$\left[k=\sum_{n=1}^{+\infty} k_n,h=\sum_{n=1}^{+\infty} h_n
 \right]=\sum_{n=1}^{+\infty}\left[ k_n,h_n
 \right]_{\cK}$$
 (one checks easily that this inner product is well defined for
 elements of $\gH_1$.)\\
 But $$\gK_{\NN}=\gH_{\NN}=\bigoplus_{n=1}^{+\infty}\cH_n=\{h=k_1\oplus k_2\oplus  ...,\; \sum_{n=1}^{\infty} ||k_n||^2_n<+\infty\}$$ as a vector space, with inner product
 $$\left[k=\sum_{n=1}^{+\infty} k_n,h=\sum_{n=1}^{+\infty} h_n
 \right]=\sum_{n=1}^{+\infty}\left[ k_n,h_n
 \right]_{\cK}$$ is a second valuable choice, distinct from the first since $e=e_1\oplus e_1\oplus  e_1...$ is in $\gH_{\NN}$ but not in $\gH_1$.\\
  In terms of
 kernels, this has the following interpretation:
Let $\cE=\RR^{\NN\times \{1,2\}}$ endowed with the topology of pointwise convergence. $\gH_1$ and $\gH_2$ are Hilbert subspaces of $\cE$, and $\gK_1$ and $\gK_2$ are distinct
Kre\u{\i}n subspaces of $\cE$ but
 with the same kernel (that is a kernel of multiplicity)
 $$K((n,.),(m,.))=\delta_{n,m}\left(
       \begin{array}{cc}
         0 &  1 \\
         1 & 0 \\
       \end{array}
     \right)$$
\end{example}
Note that by corollary 2 p.253 in \cite{Schwartz64}, $\cE$ being the dual of a barreled nuclear space, this kernel is automatically of
multiplicity.\\

The idea is then to work on fundamental Kre\u{\i}n-Hilbert pairs.\\

 Let $T$ be a locally compact space endowed with a measure
$m$. Let $t\to (\cK_t,\cH_t)$ an application from $T$ to $KH((\cF,\cE))$, $(K_t,H_t)$ be the kernels of $(\cK_t,\cH_t)$. We say that the family $\{(\cK_t,\cH_t),\; t\in T\}$ is pseudo
$m$-integrable if for all $\phi\in \cF$, the function $t\to (\phi, H_t(\phi))_{(\cF,\cE)}$ is integrable with respect to $m$. The
integral will be constructed as follows:\\

Let \begin{equation}\Pi\cH=\left\{\left\{h_t\in\cH_t,\,t\in T\right\}, \; \int_{T} ||h_t||^2_{\cH_t}<+\infty\right\}/\cR\end{equation} where $\cR$ is the equivalence relation of $m$
almost sure equality, with norm
\begin{equation}||\{h_t\in\cH_t,\,t\in T\}||^2_{\Pi\cH}= \int_{T}
||h_t||^2_{\cH_t}\end{equation} This norm makes  $\Pi\cH$ a Hilbert space. Let $\Pi\cK=\Pi\cH$ as a vector space and endow it with the indefinite form
\begin{equation}\left[\left\{k_t\in\cK_t,\,t\in T\right\},\left\{k^{'}_t\in\cK_t,\,t\in T\right\}\right]_{\Pi\cK}= \int_{T} \left[ k_t,k^{'}_t\right]_{\cK_t}dm(t)\end{equation}
Then $\Pi\cK,[.,.]_{\Pi\cK}$ is a Kre\u{\i}n space.\\

From \cite{Schwartz64} there exists a continuous linear application $\Phi$ from $\Pi\cH$ to $\cF^*$ (algebraic dual of $\cF$ endowed with the topology $\sigma(\cF^*, \cF)$) defined by:
\begin{equation}
\Phi(\{h_t\in\cH_t,\,t\in T\})=\int_T h_t dm(t)
\end{equation}
where the second member is understood as the weak integral of a scalarly integrable function. From the general theory of Hilbert subspaces \cite{Schwartz64} and subdualities
\cite{Mary2005} (and previous section) the image of the Hilbert space $\Pi\cH$ is the Hilbert subspace $\int_T \cH_t dm(t)$ of $\cF^*$ with kernel $H=\int_T H_t dm(t)\in
\cL^{+}(\overline{\cF},\cF^*)$, and the image of $\Pi\cK$ is a self-subduality (pseudo-Kre\u{\i}n subspace in the terminology of \cite{Hofmann02}) $\int_T \cK_t dm(t)$
of $\cF^*$ with kernel $H=\int_T K_t dm(t)$. Remark that $\int_T \cK_t dm(t)\subset \int_T \cH_t dm(t)$ but is not equal in general.\\

If the space $\int_T \cK_t dm(t)$ is a Kre\u{\i}n space we say that
the family is $m$-integrable. %(for instance if it is complete for the
%topology induced by $\int_T \cH_t dm(t)$).
Remark that this is equivalent with saying that $\left(\int_T \cK_t dm(t),\int_T \cH_t dm(t)\right)$ is a Kre\u{\i}n-Hilbert pair. If $\Phi$ is one-to-one, then the family is actually
$m$-integrable ($\left(\int_T \cK_t dm(t),\int_T \cH_t dm(t)\right)$ is a fundamental Kre\u{\i}n-Hilbert pair) and in this case we say that the integral is direct.

Under these assumption, the norm in $\int^{\oplus}_T \cH_t dm(t)$ is just
\begin{equation}
||\int^{\oplus}_T h_t dm(t)||^2_{\int^{\oplus}_T \cH_t dm(t)}= \int_{T} ||h_t||^2_{\cH_t}dm(t)
\end{equation}
and the indefinite inner product in $\int_T \cK_t dm(t)$ is
\begin{equation}
\left[\int^{\oplus}_T k_t dm(t),\int^{\oplus}_T k^{'}_t dm(t)\right]_{\int^{\oplus}_T \cK_t dm(t)}= \int_{T} \left[k_t,k^{'}_t\right]_{\cK_t}dm(t)
\end{equation}

\section{Group representations in Kre\u{\i}n subspaces}

\subsection{Unitary representations and Hilbert subspaces}
The study of group representation in Hilbert subspaces is one of the many approaches to harmonic analysis. The usual setting is the study representations of a Lie group $G$ on Hilbert
subspaces (for instance $L^2(\mu)$ for an invariant measure on $G$) of the space of distributions on $X=G$ (or $X$ an homogeneous space \cite{Capelle96}, \cite{Thomas05}). One can also
study Hilbert subspaces of other locally convex spaces, such as spaces of
holomorphic functions \cite{Faraut99}, \cite{Faraut03}. \\
We refer to the works \cite{Carey78}, \cite{Kunze67}, \cite{Bekka03} for fundamentals theorems on irreducibility in reproducing kernel Hilbert spaces or Hilbert subspaces.
%and \cite{Faraut99}, \cite{Faraut03},
%\cite{Thomas99}, \cite{Capelle96} for
% fondamentaux -irreducibility (De
%la harpe, Kunze, Carey) -Gelfand Pairs, multiplicity free (Faraut,
%Thomas, capelle) -Zonal spherical functions

\subsection{Extension to Kre\u{\i}n subspaces}
In this section $T:\cE\rightarrow \cE$ is a weakly continuous
operator, and %$\tau:g\in G\rightarrow \tau_g$ is a representation of
$G$ is a group of weakly continuous endomorphisms of $\cE$ (in this section we identify a group $G$ with is image under a given representation $\tau$ in $\cL(\cE)$).
%Continuity of the representation is not assumed, and this topic will
%be addressed in a forthcoming paper.

\begin{definition}
A Kre\u{\i}n subspace $\cK$ (with kernel $K$) is invariant under $T$ if $T(\cK)\subset \cK$, and $T$ is a unitary operator ($T^{[*]}T=TT^{[*]}=I_{\cK}$).
\end{definition}

\begin{example}
Let $\theta\in \RR$, and define  $T_{\theta}:\RR^2\rightarrow \RR^2$ by
$$T_{\theta}=\left(
       \begin{array}{cc}
         i\sinh(\theta) &  \cosh(\theta) \\
         \cosh(\theta) & i\sinh(\theta)
       \end{array}
     \right)$$
Then the Kre\u{\i}n space $\cK$ of example \ref{ex 1} is invariant with respect to $T_{\theta}$.
\end{example}

\begin{proposition}\label{invar1}
A Kre\u{\i}n subspace $\cK$ (with kernel $K$) is invariant under $T$ if $T(\cK)\subset \cK$, $T K T^*=K$ and $T(\cK)$ is a Kre\u{\i}n space.
\end{proposition}

\begin{proof}
Let $T(\cK)$ be the inner product space, image of the Kre\u{\i}n space $\cK$. Since it is Kre\u{\i}n space, by theorem \ref{image} it is a Kre\u{\i}n subspace of $\cE$ with kernel is
$T K T*=K$. But $T(\cK)\subset \cK$, and by proposition 39 in \cite{Schwartz64} the two Kre\u{\i}n spaces coincide, and $T$ acts as a unitary operator on $\cK$.
\end{proof}

Of course, if $\cK$ is finite-dimensional then the equality $T K T^*=K$ is sufficient.

\begin{definition}
A Kre\u{\i}n subspace $\cK$ (with kernel $K$) is invariant under
%(the action $\tau:G\to \cL(\cK)$ of a group)
$G$ if $\forall g\in G,\;g(\cK)\subset \cK$, $g K g^{*}=K$.
\end{definition}

\begin{remark}\label{leK}
Since $G$ is a group, then $g(\cK)$ is always a Kre\u{\i}n space.
\end{remark}

\begin{theorem}
A Kre\u{\i}n subspace $\cK$ (with kernel $K$ and minimal bound $H$) is invariant under a group $G$ if and only if
\begin{enumerate}
\item $\forall g\in G,\;g
K g*=K$
\item $\forall g\in G,\;\exists \alpha_g,\;g
H g*\leq \alpha_g H$
\end{enumerate}
\end{theorem}

\begin{proof}
Let $g\in G$. By proposition \ref{invar1} and remark \ref{leK}, $\cK$ is invariant under $g$ if $g(\cK)\subset \cK$ or equivalently, if $g(\cH)\subset \cH$. But $g(\cH)$ is a Hilbert
subspace of $\cE$ with kernel $gHg^*$, and by proposition 15 in \cite{Schwartz64}, the inclusion holds if and only if exits $c>0,\; gHg^*\leq c H$.
\end{proof}

%Principaux theorem a reprendre dans le cas des SD. Probleme des
%noyaux minimaux ! Travailler sur les paires de noyaux.

\begin{definition}
A invariant Kre\u{\i}n subspace $\cK$ %$(\cK,[.,.])$
is called (topologically) irreducible if it admits no invariant closed subspace
apart from $0$ and $\cK$.\\
It is called indecomposable if any decomposition as a direct sum of two invariant closed subspaces involves the trivial subspaces.\end{definition}

\begin{definition}
A invariant Kre\u{\i}n subspace $\cK$ %$(\cK,[.,.])$
is called regularly irreducible if it admits no invariant regular subspace (closed and a Kre\u{\i}n space with the induced sesquilinear
form) apart from $0$ and $\cK$.\\
It is called regularly indecomposable if any decomposition as a direct sum of two invariant regular subspaces involves the trivial subspaces.\end{definition}

%\begin{definition}
%A invariant Kre\u{\i}n subspace $\cK$ %$(\cK,[.,.])$
%is called contractively irreducible if it admits no invariant
%contractive subspace (Kre\u{\i}n space contractively
%embedded in $\cK$) apart from $0$ and $\cK$ ($\lambda \cK$ if the inner product is allowed to be definite) .\\
%It is called contractively indecomposable if any decomposition as a
%complementary sum \cite{DeBranges88} of two invariant contractive
%subspaces involves the trivial subspaces.\end{definition}

Finally, following Kissin\cite{Kissin96}, we say that a
representation %$\tau$ of $G$
on $\cK$ is non-degenerate if $\cK$ has no neutral invariant subspace.

From \cite{Schwartz64Lie} we have:
\begin{lemma}
%For a Hilbert space, the three notions of irreducibility coincide.
For a Hilbert space, the two notions of irreducibility coincide.
\end{lemma}

\begin{theorem}
For Krein spaces: \begin{itemize} \item Irreducibility $\Rightarrow$ indecomposability, but the converse is not true.
%Regular (contractive) irreducibility $\iff$ regular (contractive)
\item Regular irreducibility $\iff$ regular indecomposability.
\end{itemize}
\end{theorem}

\begin{proof}
The implication is  straightforward.\\
For the converse, suppose that $\cK$ is indecomposable, and let $\cH$ be a regular subspace of $\cK$. Then $\cH^{\perp}$ is a regular subspace and $\cH\oplus \cH^{\perp}=\cK$. By
indecomposability, $\cH=0$ or $\cH=\cK$.%\\
%For contractions, we use the complementation theory of de
%Branges\cite{DeBranges88}.
\end{proof}

\subsection{Schur's lemma in Kre\u{\i}n spaces}
A main tool in harmonic analysis is Schur's lemma, that asserts that for unitary representations in Hilbert spaces irreducibility is equivalent with operator irreducibility. The aim of
this section is to study the link between the different notions in the Kre\u{\i}n space setting. Note that a strange phenomenon will occur, for in the algebra of bounded operators on a
Kre\u{\i}n space the Gelfand-Naimark property is not valid \cite{Mary2008MPK}, and there exists
self-adjoint nilpotent operators.\\

Also, the set of self-adjoint operators in a Kre\u{\i}n space is too large for a good spectral theory. Hence it is classical to study definitizable operators \cite{Iokhvidov82},
\cite{Langer82}. An operator $A$ is definitizable if there exists a polynomial $p$ such that $p(A)$ is a positive operator. Hence we start with the study of positive operators in the
commutant of $G$.

\begin{lemma}\label{lemmaPosReg}
Let $\cK$ be a Kre\u{\i}n space, regularly irreducible under a the action
of a group of unitary operators $G$.\\
Let $A$ be a bounded positive operator and suppose that
$$\forall T\in G,,\; TA=AT$$ Then either
$\cK$ is a Hilbert space and $A=\lambda I,\; \lambda \geq 0$ or $\cK$ is not definite and $A=N; \; N \textit{nilpotent and
positive}, N^2=0$.\\

If moreover the representation is non-degenerate then $A=0$.
\end{lemma}

\begin{proof}
If $\cK$ is a Hilbert space then it is Schur's lemma.\\ Suppose now that $\cK$ is not definite. Since $A$ is positive we have a spectral decomposition \cite{Langer82}. We note
$\cK_{\rho}$ the corresponding spectral subspaces. These are Kre\u{\i}n spaces and there orthogonal sum is $\cK$. Moreover, these spaces are invariant since spectral projections are
in the double commutant of $A$.\\ %(see Theorem 2.1 ([33] in indefinite
%Sturm-Liouville... ou Theorem A7 p147 in
%Methoden der Hartree-Fock-Bogoliubov Theorie).\\
It follows that all of the spectral subspaces are $\{0\}$ except one that is $\cK$. By spectral decomposition, $A-\rho I$ is nilpotent, $A=\rho I+N$. Suppose $\rho\neq 0$. Then
$B=I+\rho^{-1}N=I+M$ is positive and invertible, $B^{-1}=I-M$ is also positive and $B+B^{-1}=2I$ is positive, which is excluded by hypothesis. It follows that $\rho=0$, and $N$ is
positive, $N^2=0$.
\end{proof}

%\begin{lemma}
%Let $\cK$ be a Kre\u{\i}n space, contractively irreducible under a
%the action
%of a group of bounded operators $G$.\\
%Let $A$ be a bounded positive operator %(see generalized projections
%of Curgus and Langer, complementations in Kre\u{\i}n spaces),
%and suppose
%that
%$$\forall T\in G,,\; TA=AT$$ Then either
%$\cK$ is a Hilbert space and $A=\lambda I,\; \lambda \geq 0$ or
%$\cK$ is not definite and $A=0$.
%\end{lemma}

%\begin{proof}
%By lemma \ref{lemmaPosReg}, either $A=\lambda I,\; \lambda \geq 0$
%and $A$ is not positive unless $\cK$ is a Hilbert space, or $A=N$ is
%nilpotent and positive. Suppose $N\neq 0$. Then it defines a non
%trivial Kre\u{\i}n space $\cN$ contractively included in $\cK$, and
%$\cK$ is not contractively irreducible.
%\end{proof}

We can now state a general result for definitizable operators in the commutant of $G$:
\begin{lemma}[Regular Schur's lemma for Kre\u{\i}n
spaces]\label{RSLemma} Let $\cK$ be a Kre\u{\i}n space, regularly irreducible under a the action
of a group of unitary operators $G$.\\
Let $A$ be a bounded definitizable (self-adjoint) operator
%(see
%generalized projections of Curgus and Langer, complementations in
%Kre\u{\i}n spaces),
and suppose that
$$\forall T\in G,\; TA=AT$$ Then
$$A=N+\lambda I; \; N \textit{nilpotent.}$$
If moreover the representation is non-degenerate then $A=\lambda I$.
\end{lemma}

\begin{proof}
Since $A$ is definitizable we have a spectral decomposition.
%We note
%$\cK_{\rho}$ the corresponding spectral subspaces. These are
%Kre\u{\i}n spaces and there orthogonal sum is $\cK$. Moreover, these
%spaces are invariant since spectral projections are in the double
%commutant of $A$ (see Theorem 2.1 ([33] in indefinite
%Sturm-Liouville... ou Theorem A7 p147 in
%Methoden der Hartree-Fock-Bogoliubov Theorie).\\
As before, by irreducibility all of the spectral subspaces are $\{0\}$ except one
that is $\cK$, and $A-\rho I$ is nilpotent.\\

%Let $\mathfrak{H}_1=\overline{Im A}=(\ker
%A)^{\perp},\;\mathfrak{H}_2=(Im A)^{\perp}=\ker A$.

%NON : If $\exists x\in \mathfrak{H}_2 \backslash \mathfrak{H}_1$,
%then $\CC.x$ is invariant under $G$ and $[x,x]\neq 0$ hence it is a
%regular subspace of $\cK$. Absurd. So
%$\mathfrak{H}_2=\mathfrak{H}_1$. NON\\

Suppose now that the representation is non-degenerate. Let $N$ be nilpotent of order $k>1$ and pose $M=N^{k-1}$. Then $M^2=0$, $\overline{Im M}$ is neutral and invariant hence $0$ or
$\cK$, absurd. Then $k=1$ and $N=0$.
%It is also the case if there
%exists no (non-trivial) invariant Hilbert subspace (seulement dans
%le cas d'un polynome a racines simples !).

\end{proof}

\begin{example}\label{ex 3}
Let $\cE=\cH\oplus \cH$ where $\cH$ is a Hilbert space and
$$G=\{\left(
                  \begin{array}{cc}
                    I & B \\
                    0 & I \\
                  \end{array}
                \right)\}$$ with $B$ anti-symmetric. Remark that $G$ is a (multiplicative) group. Suppose also that an invertible
                antisymmetric operator exists.

            Pose $J=\left(\begin{array}{cc}
                                                           0 & 1 \\
                                                           1 & 0 \\
                                                         \end{array}
                                                       \right)$.
             Then $G$ defines J-unitary operators, $Q=\left(\begin{array}{cc}
                                                              0 & 1 \\
                                                              0 & 0 \\
                                                            \end{array}
                                                          \right)$
             is nilpotent and commutes with $G$. Remark that %$\cH$ is an invariant unitary representation
             %of $G$,
             $\cM=(\cH,0)$ is invariant under $G$, neutral in $\cK=\cH\oplus \cH$ and
             $\cK$ is regularly irreducible, since only trivial (J-self-adjoint) projection
             commute with $G$ by the existence of the invertible and antisymmetric
             operator.
             %Remark that the representation is degenerate, $\cM=\CC e_1$
             %is invariant and neutral.
\end{example}
%\begin{theorem}
%Regular irreducibilty $\iff$ pair $(K, L)$ minimal.
%\end{theorem}

%\begin{proof}
%For the implication, we make use of the functional calculus for
%normal operators. (see coherent states, or Knapp representation
%theory of semisimple groups).

%The converse is straghtforward using orthogonal projection.
%\end{proof}

%NON NE MARCHE PAS. IL FAUT $JA=AJ$ ou $A=ii^*$. Soit $H_+\subset
%K_+$ et $H_-\subset K_-$

%\begin{corollary}[Commutant's regular Schur's lemma]
%
%\end{corollary}
\subsection{``Fundamental'' representations}

It may happen that the %(representation $\tau$ of the)
group $G$ carries a fundamental symmetry of $\cK$: \begin{equation} \exists
%g\in G,\; \tau_g=\tau_g^{[*]}=(\tau_g)^{-1} \text{ and } \tau_g
g\in G,\; \g=g^{[*]}=g^{-1} \text{ and } g \text{ positive.}
\end{equation}
In this case, we say that the representation is fundamental. The interest of fundamental representations lies in the following lemma and theorem:

\begin{lemma}[Regular Schur's lemma for fundamental representation]
Let $\cK$ be a Kre\u{\i}n space, regularly irreducible under the ``fundamental'' action
of a group of unitary operators $G$.\\
Let $A$ be a bounded operator
%(see
%generalized projections of Curgus and Langer, complementations in
%Kre\u{\i}n spaces),
and suppose that
$$\forall T\in G,\; TA=AT$$ Then
$$A=\lambda I$$
\end{lemma}

%In other words, for fundamental representations regular
%irreducibility is equivalent with commutant irreducibilty.

\begin{proof}
Let $J$ be a fundamental symmetry associated to $G$. Pose $B=\frac{A+A^{[*]}}{2}$, $C=i\frac{A-A^{[*]}}{2}$. Since $G$ acts by unitary operators, $A^[*]$ also commutes with $G$ and $B,
C$ are self-adjoint and commute with $G$. But they also commute with a fundamental symmetry, hence they are self-adjoint in the Hilbert space sense, admit a spectral function and by
irreducibility their spectrum reduces to a single number. Finally $A=\lambda I$.
\end{proof}

\begin{theorem}
For fundamental representations, regular irreducibility implies (topological) irreducibility.
\end{theorem}

\begin{proof}
Let $\cL$ be a closed subspace invariant with respect to the fundamental symmetry $J\in G$. Then $J\cL=\cL$, and it follows that $\cL^{[\perp]}=(J\cL)^{\perp}=\cL^{\perp}$. But
$\cL\oplus \cL^{\perp}=\cK$ and $\cL$ is regular, hence trivial.
%but $J$ is positive hence $x=0$ and $\cM=\{0\}$.
\end{proof}
%Remark that a fundamental representation is non-degenerate:
%In fact, we could have also used the non-degeneracy property:
%\begin{proposition}
%A fundamental representation is non-degenerate.
%\end{proposition}

%\begin{proof}
%Let $\cM$ be a neutral subspace invariant with respect to the
%fundamental symmetry $J=\tau_g$. Then $\forall x\in \cM,\; [x,Jx]=0$
%but $J$ is positive hence $x=0$ and $\cM=\{0\}$.
%\end{proof}

As a consequence:
\begin{corollary}
For fundamental representations, the three following notions coincide:
\begin{enumerate}
\item Regular irreducibility; \item Operator
irreducibility; \item Topological irreducibility.
\end{enumerate}
\end{corollary}

Note that a fundamental representation is then obviously non-degenerate.

%\begin{proof}
%\begin{enumerate} \item For the implication this is Schur's lemma. For the converse, note
%that an invariant regular subspace is defined by a projection in the
%commutant of $G$.
%\item Suppose $\cK_0$ is contractively included in $\cK$ and
%invariant. It follows that tts kernel $P$ on $\cK$ commutes with any
%$T\in G$. By Schur's lemma, $P=\lambda Id$ and finally $P=0$ or
%$P=Id$ and $\cK_0$ is trivial.
%\end{enumerate}
%\end{proof}

\subsection{Reproducing kernel Kre\u{\i}n space and irreducible representations}
In this section we suppose that $\cE=\CC^{X}$, where $X$ is a set, and $G$ is a group acting transitively on $X$ (for instance $X$ is an homogeneous space). There is a canonical action
of $G$ on $\cE$ defined by: \begin{equation} \forall g\in G,\; \tau_g(f)(x)=f(g^{-1}(x))
\end{equation}
Let $\cK$ be a Kre\u{\i}n subspace of $\cE$ (we call such a subspace a reproducing kernel Kre\u{\i}n space) invariant with respect to
$\tau$. \\

Note that for such a subspace \cite{Schwartz64}, \cite{Mary2003}, \cite{Mary2005} we can identify its (unique) kernel with a reproducing function $K(.,.)$ on $X^2$ that verifies:
\begin{equation}
\forall (x,y)\in X^2,\; \left[K(x,.),K(y,.)\right]_{\cK}=K(x,y)
\end{equation}
or equivalently (this is equation \ref{eqK})
\begin{equation}
\forall x\in X,\;\forall k\in \cK,\; \left[K(x,.),k\right]_{\cK}=k(x)
\end{equation}

For any $\omega\in X$, we define its isotropy group $\varpi=\{g\in G,\; g\omega=\omega\}$. For such a subgroup we define the subspace $\cK^{\varpi}=\{k\in \cK,\; \forall \varrho\in
\varpi,\; \tau_{\varrho}(k)=k\}$ of $\varpi$-invariant functions.
\begin{theorem}\label{th_class1}
$\,$\\ \vspace{-1cm}
\begin{enumerate}
\item if $K$ is the reproducing kernel function of $\cK$,
\begin{equation}\label{equ_Kinvar1}
\forall g\in G,\;\forall (x,y)\in X^2,\; K(gx,gy)=K(x,y)
\end{equation}
\item if $\cK\neq \{0\}$ then $\cK^{\varpi}\neq \{0\}$.
\item if $\dim(\cK^{\varpi})=1$ then the representation is
regularly irreducible.
\end{enumerate}
\end{theorem}

\begin{proof}
\begin{enumerate}
\item
Fix $g\in G$ and define $R(x,y)=K(gx,gy)$. Since the representation is unitary,
\begin{eqnarray*}
\left[R(x,.),k\right]_{\cK}&=&\left[K(gx,g.),k\right]_{\cK}=\left[\tau_g(K(gx,g.)),\tau_g(k)\right]_{\cK}\\
&=&\left[K(gx,.),k(g^{-1}.)\right]_{\cK} =k(g^{-1}gx)\\
&=&k(x)
\end{eqnarray*}
and we conclude by unicity of the kernel.
\item Since $\cK\neq \{0\}$, the function $K(.,.)$ is not identically zero and exists $x,y$ in $X^2$, $K(x,y)\neq 0$. Since $G$ acts transitively on
$X$, exists $g\in G,\; gx=\omega$ and by equation \ref{equ_Kinvar1} $K(\omega,gy)\neq 0$. But still by equation \ref{equ_Kinvar1} and the definition of $\varpi$ the function
$k(.)=K(\omega,.)$ is in $\cK^{\varpi}$, and $\cK^{\varpi}\neq \{0\}$.
\item Let $\cK_0$ be a regular
subspace of $\cK$. Then it is a Kre\u{\i}n space continuously included in $\CC^{X}$ hence it admits a reproducing kernel function $K_0$. But
$K_0$ is then $G$-invariant and %the function
$K_0(\omega,.)\in \cK^{\varpi}$. Since $\dim(\cK^{\varpi})=1$ $K_0(\omega,.)$ is proportional to $K(\omega,.)$ and by transitivity of $G$, $K_0$ is proportional to $K$:
\begin{equation}
\exists c\in \CC,\; \forall x,y\in X^2,\; K_0(x,y)=c K(x,y)
\end{equation}
Now let $k_0\in \cK_0,\; k_0(x)\neq 0$.
\begin{eqnarray*}
k_0(x)&=&\left[ K_0(x,.),k_0(.)\right]_{\cK}\\
&=&\overline{c}\left[ K(x,.),k_0(.)\right]_{\cK}\\
&=&\overline{c}k_0(x)
\end{eqnarray*}
and $c=1$. The two Kre\u{\i}n subspaces have the same kernel, one is included in tho other hence by proposition 39 in \cite{Schwartz64} they coincide as Kre\u{\i}n spaces.
\end{enumerate}
\end{proof}

\begin{example}[Homogeneous polynomial representations of the Lorentz
group]\label{ex HP} Let $X=\RR^3$ and $G=SO(1,2)$ be the associated Lorentz group. The action of $G$ on $X$ is transitive, and if $J=\left(
                                                 \begin{array}{ccc}
                                                   1 & 0 & 0 \\
                                                   0 & -1 & 0 \\
                                                   0 & 0 & -1 \\
                                                 \end{array}
                                               \right)$
is the Minkowski metric operator, then the reproducing kernel $K(x,y)=(\langle Jx,y\rangle_{\RR^3})^n$ is invariant under $G$ and it defines a finite-dimensional Kre\u{\i}n space of
homogeneous
polynomials of degree $n$.\\
Let $\omega=e_3$. Then its isotropy subgroup $\varpi$ contains the Lorentz boosts $$T_{\theta}=\left(
                                                 \begin{array}{ccc}
                                                   \cosh(\theta) & \sinh(\theta) & 0 \\
                                                   \sinh(\theta) & \cosh(\theta) & 0 \\
                                                   0 & 0 & 1 \\
                                                 \end{array}
                                               \right)$$
It follows the the space of $\varpi$ invariant functions is at most
one-dimensional:\\
Let $k$ be an invariant homogeneous polynomial of degree $n$. Then it is of the form \begin{equation} k(x)=\sum_{i\in I}\alpha_i (\langle, Jy_i, x\rangle)^n
\end{equation}
and its invariance under the Lorentz boosts $T_{\theta}$ implies
$$\forall i\in I, \; T_{\theta}y_i=y_i$$
Finally $y_i\propto(0,0,1)$ %an polynomial of degree $n$ invariant
%under all such Lorentz boosts
and $k$ is of the form $k(x_1,x_2,x_3)=a (x_3)^n$.\\
But $\cK^{\varpi}$ is also at least of degree one by theorem \ref{th_class1} $(2)$, and by $(3)$ $\cK$ is regularly irreducible. Note that the representation is actually fundamental
($J\in G$), and $\cK$ is topologically irreducible.
\end{example}

\section{Integral decomposition in convex cones}
In this section, any convex cone $\Delta$ will induce its proper order $\leq_{\Delta}$: \begin{equation}(d,d')\in \Delta^2,\; d'\leq_\Delta d\iff d-d'\in \Delta\end{equation} Two
elements $d,d'$ will be called $\Delta$-independent if $h\leq_{\Delta} d$, $h\leq_{\Delta} d'$ implies $h=0$ and we note $d \amalg_{\Delta} d'$. This is the relation we definition we
used for the cone of positive kernels.

\subsection{Integral representation property for closed convex cones}
%We say that a closed convex cone has the I.R.P. if

First we recall the definition of the integral representation property (I.R.P.), and the main theorem of integral representation in conuclear cones due to Thomas \cite{Thomas94}. In
the following $ext(\Gamma)$ denotes the set of extreme rays of any closed convex cone $\Gamma$, and $\cM^+(X)$ the set of positive radon measures on the topological set $X$.

\begin{definition}
A closed convex cone $\Gamma$ has the I.R.P. if
\begin{enumerate}
\item for every closed convex subbcone $\Gamma_1\subset \Gamma$, the
map $r:\cM^+(ext(\Gamma_1))\rightarrow \Gamma_1$ is onto;
\item the map $r:\cM^+(ext(\Gamma_1))\rightarrow \Gamma_1$ is
bijective if and only if $\Gamma_1$ is a lattice (with respect to its proper order).
\end{enumerate}
\end{definition}

\begin{theorem}
Let $F$ be a weakly complete conuclear space. Then any salient and closed convex cone $\Gamma\subset F$ has the I.R.P.
%\begin{enumerate}
%\item $\Gamma$ is the closed convex hull of its extreme rays.
%\item If $t\to (e_t)$, $T\to ext(\Gamma)$ is an admissible
%parametrization of the extreme rays then :\begin{enumerate}
%\item for every $f\in Gamma$ there exists a Radon
%measure $m$ on $T$ such that
%\begin{equation}f=\int_T (e_t dm(t)\end{equation}
%\item the measure $m$ is uniquely determined by $f$ if and only if the face
%$\Gamma(f)$ generated by $f$ is a lattice.
%\end{enumerate}
%\end{enumerate}
\end{theorem}

For the rest of this section we suppose that $F=\mathfrak{F}\times \mathfrak{F}$ is a weakly complete conuclear space, $\gC$ is a salient closed convex cone of $\gF$ and $\gV=\gC-\gC$
is the vector space generated by $\gC$. It follows that $\digamma=\gV\times \gC$ is a salient closed convex cone of $F$ that has the I.R.P. (for our application, $\gF$ will be the
space of kernels and $\gC$ the cone
of positive kernels) and every closed convex subcone (for instance the pairs of invariant kernels) will have the I.R.P.\\

%Finally we consider three more cones: a closed convex subcone
%$\Gamma\subset \gV\times \gC$, the closed convex cone of dominated
%pairs $D=\{(v,c)\in \gV\times \gC,\; -c\leq_{\gC} v \leq_{\gC} c\}$
%\footnote{the order is the order induced by the convex cone $\gC$, $v\leq c\iff c-v\in \gC$}
%and finally the intersection of these two cones
%$\Gamma_D=\Gamma\cap D$.\\

%We make the following assumptions:
%\begin{itemize}
%\item[(1)] $\Gamma$ is of the form $\Gamma=\Gamma_1\times \Gamma_2$
%where $\Gamma_2$ is a convex cone hereditary for the order induced
%by $\gC$ and $\Gamma_1$ is a vector space included in $\Gamma_2-
%\Gamma_2$.

%\item[(2)] If $(c,d)\in \Gamma_2^2$ are  $\Gamma_2$-independent, then
%they are $\gC$-independent.

%\item[(3)]
%\end{itemize}

Let $D=\{(v,c)\in \gV\times \gC,\; -c\leq_{\gC} v \leq_{\gC} c\}$ be the closed convex cone of dominated pairs. We now define the set of minimal pair in this setting:
\begin{definition}
The pair $(v,c)\in D$ is $\gC$-minimal (or a minimal pair) if $c-v$ and $c+v$ are independent for the order induced by $\gC$ ($c-v\amalg_{\gC} c+v$), equivalently if any $h\in \gC$
verifying $h\leq_{\gC} c-v,\; h\leq_{\gC} c+v$ is zero.
\end{definition}

\begin{lemma}\label{mini}
Let $(v,c)\in D$ be a minimal pair. If $(w,d)\leq_{D} (v,c)$, then $(w,d)$ is $\gC$-minimal.
\end{lemma}

\begin{proof}
Let $(u,p)=(v-w,c-d)\in D$. Let $h\leq_{\gC} d-w,\; h\leq_{\gC}
d+w$. Then %\begin{equation}-(d-h)-p=h-d-p\leq_{\gC} w+u \leq_{\gC} (d-h)+p\end{equation}
\begin{equation}h-d-p\leq_{\gC} w+u \leq_{\gC} (d-h)+p\end{equation}
and since $w+u=v,\; d+p=c$, $h-c\leq_{\gC} v \leq_{\gC} c-h$ which gives \begin{equation}h\leq_{\gC} v+c,\; h\leq_{\gC} c-v\end{equation} and $h=0$ since $(v,c)$ is minimal.
\end{proof}

\begin{lemma}\label{lh}
Let $(v,c)\in D$ not be minimal. Then there exists $h\in \gC-\{0\}$ and $(v,d)\in D$, $(v,c)=(v,d)+(0,h)$
\end{lemma}

Let $\Gamma$ be a closed convex subcone of $\digamma$ and $(v,c)\in \Gamma_D=\Gamma\cap D$ be a minimal pair. By the previous theorems if $t\to (e_t,f_t)$, $T\to ext(\Gamma_D)$ is an
admissible parametrization of the extreme rays then there exists a Radon measure $m$ on $T$ (unique if the face $\Gamma_D\left(\left(v,c\right)\right)$ is a lattice) such that
\begin{equation}\label{eq_int_kern}(v,c)=\int_T (e_t,f_t) dm(t)\end{equation}

\begin{lemma}\label{asM}
%If $\gF'$ is separable then the set of $(t,t')\in T^2$ such that
The set of $(t,t')\in T^2$ such that $(e_t+e_{t'}, f_t+f_{t'})$ is not a minimal pair is of $m$ measure $0$.
\end{lemma}

\begin{proof}
Suppose $m$ is of mass one. We have \begin{equation}(v,c)=\int_T (e_t,f_t) dm(t)=\int_{T^2} (e_t+e_{t'},f_t+f_{t'}) d(m\otimes m) (t,t')\end{equation} From the construction in
\cite{Thomas94}, $m$ is concentrated on a compact and metrizable (hence separable) set of
$B$ of $\gF'$, and %any metrizable and compact set being separable,
exists $\{\varphi_n, n\in \NN\}$ a dense family of $B$. Let
\begin{equation}N=\{(t,t')\in T^2,\; (e_t+e_{t'}, f_t+f_{t'}) \text{ not
minimal} \}\end{equation} Then by lemma \ref{lh}
\begin{equation}\forall (t,t')\in N,\; \exists
h_{t,t'}>0,\;(e_t+e_{t'}, f_t+f_{t'})=\gamma_{t,t'}+(0,\h_{t,t'})\end{equation} Define
\begin{equation} N_n=\{(t,t')\in N,\; \varphi_n(h_{t,t'})>0\}\end{equation}
%\begin{equation}N_n=\{(t,t')\in T^2,\; (e_t+e_{t'}, f_t+f_{t'})=\gamma_{t,t'}+(0,\h_{t,t'}),\; \gamma_t \text{ minimal and }\varphi_n(h_{t,t'})>0\}\end{equation}
Then $N=\bigcup_{n\in \NN} N_n$ by the Hahn-Banach theorem.
\begin{equation}
\int_{N_n} (e_t+e_{t'},f_t+f_{t'}) d(m\otimes m) (t,t')=\int_{N_n} \gamma_{t,t'} d(m\otimes m) (t,t')+\int_{N_n} (0,h_{t,t'}) d(m\otimes m) (t,t') \end{equation} and by lemma
\ref{mini} $\int_{N_n} (0,h_t) d(m\otimes m) (t,t')=(0,h)$ is a minimal pair hence $\int_{N_n} h_{t,t'} d(m\otimes m) (t,t')=0$. It follows that $\int_{N_n} \varphi_n(h_{t,t'})
d(m\otimes m) (t,t')=0$ which implies that $m(N_n)=0$ since $\varphi_n(h_{t,t'})>0$ (This is a classical application of the monotone convergence theorem). Finally by
$\sigma-$subadditivity $m(N)=0$.
\end{proof}

\section{Applications to invariant kernels and Kre\u{\i}n subspaces}
 As for the direct integral of Kre\u{\i}n subspaces, where the
Kre\u{\i}n space structure only was not sufficient, the kernel alone of a Kre\u{\i}n subspace is not sufficient to have a minimal decomposition. This is due to the vector space
structure of the set of hermitian kernels, hence the fact the order intervals are not bounded. To get an integral decomposition, we work on Kre\u{\i}n-Hilbert pairs and their kernels.
.
\subsection{Integral decomposition of invariant minimal Kre\u{\i}n-Hilbert pairs}
Let $\gC=\bL^{+}(\cF,\cE)$ be the cone of positive kernels in $\cL{(\cF,\cE)}$ and $\gV=\gC-\gC=\bL^{b}(\cF,\cE)$ be the vector space generated by $\gC$, vector space of bounded
hermitian kernels. Suppose that the space $\cL{(\cF,\cE)}^2$ is weakly complete and conuclear. Then the cone $\digamma=\gV\times \gC$ is a salient closed convex cone, and it has the
I.R.P. (this will for notably be the case if $\cE$ is the space of distributions on a Lie group $G$ \cite{Thomas94}). Let $\cT$ be a group of weakly continuous operators on $\cE$, and
define the following convex cones:

\begin{equation}\label {UU} \UU=\{K\in \gV,\; \forall T\in \cT,\; TKT^*=K \}
\end{equation} and for $\Lambda=\{\lambda_T,\; T \in \cT\}$ a family
of positive numbers indexed by $\cT$
\begin{equation} \II_{\Lambda}=\{H\in \gC,\; \forall T\in \cT,\; THT^*\leq \lambda_{T} H \}\end{equation} Finally let
$\Gamma_D(\Lambda)=(\UU\times \II_{\Lambda})\bigcap D$ where $D$ is the cone of dominated pairs.

\begin{lemma}
$\Gamma_D(\Lambda)$ is a closed convex subcone of $\digamma$.
\end{lemma}

\begin{proof}
The convexity is straightforward, and the closedness follows from
the weak continuity of $T\subset \cT$: \\
For all $\gamma\in \RR$, the operators $\widehat{\gamma T}:\cL{(\cF,\cE)}\to \cL{(\cF,\cE)}$ defined by $\widehat{\gamma T}(v)=\gamma.v-TvT^*$ are weakly continuous, $\UU=\bigcap_{T\in
\cT}\widehat{1 T}^{-1}\{0\}$ and $\II_{\Lambda}=\bigcap_{T\in \cT}\widehat{\lambda_T T}^{-1}\{\gC\}$
\end{proof}

Let $(\cK,\cH)$ be an invariant minimal Kre\u{\i}n-Hilbert pair of subspaces of $\cE$ with kernels $(K,H)$. Using the results of section 2 $\exists \Lambda, (K,H)\in
\Gamma_D(\Lambda)$. By the previous theorem if $t\to (K_t,H_t)$, $T\to ext(\Gamma_D(\Lambda))$ is an admissible parametrization of the extreme rays then there exists a Radon measure
$m$ on $T$ (unique if the face $\Gamma_D(\lambda)\left(\left(K,H\right)\right)$ is a lattice) such that
\begin{equation}\label{eq_int_kern}(K,H)=\int_T (K_t,H_t) dm(t)\end{equation}

Define $(\cK_t,\cH_t)$ the associated family of invariant Kre\u{\i}n-Hilbert pairs of subspaces of $\cE$ with kernels $(K_t,H_t)$. Then by equation \ref{eq_int_kern}, this family  is
pseudo $m$-integrable and we can define their integral $\int_T (\cK_t,\cH_t) dm(t)$.

\begin{proposition}
If the family is $m$-integrable then \begin{equation} \int_T (\cK_t,\cH_t) dm(t)=(\cK,\cH)\end{equation}
\end{proposition}

\begin{proof}
the space $\int_T \cK_t dm(t)\subset \int_T \cH_t dm(t)=\cH$, but $\cH=\cK$ as subspaces and since $\int_T \cK_t dm(t)$ and $\cK$ are Kre\u{\i}n subspaces, then they coincide
(proposition 39 p~246 in \cite{Schwartz64}).
\end{proof}

In fact we have more :

\begin{theorem}\label{direct}
The integral $\int_T (\cK_t,\cH_t) dm(t)$ is a direct integral and
\begin{equation}\label{eq direct integral of spaces} (\cK,\cH)=\int^{\oplus}_T
(\cK_t,\cH_t) dm(t)\end{equation}
\end{theorem}

% in the following sense: if $A,B$
%are two disjoint compact subsets of $T$, then $\int_A f_t dm(t)$ and
%$\int_B f_t dm(t)$ are independent.

%\begin{lemma} FAUX
%Let $(v,c)\in \Gamma_D$ be extremal. Then it is a minimal pair.
%\end{lemma}

%\begin{proof}
%Let $h\leq_{\gC} c-v, c+v$. Then $h\in Gamma_2$ by heredity of
%$\Gamma_2$ hence $(0,h)\in \Gamma$ and $(0,h)\leq_{\Gamma_D} (v,c)$,
%hence $h=0$ by extremality.
%\end{proof}

To prove this theorem we need the following lemma

\begin{lemma}\label{interM}
Let $(V,C)$, $(W,D)$  and $(U,P)$ be three minimal pairs in $\Gamma_D(\Lambda)$ such that $(V,C)=(W,D)+(U,P)$. Suppose moreover that $(W,D)$, $(U,P)$ are extremal. Then either they
belong to the same extreme ray or $D$ and $P$ are $\gC$-independent.
\end{lemma}

\begin{proof}
\begin{eqnarray}
& & C=D+P,\; v=W+U\\
&\Rightarrow& (C-v)=(D-W)+(P-U),\; (C+v)=(D+W)+(P+U)\\
&\Rightarrow& (D-W)+(P-U)\amalg_{\gC} (D+W)+(P+U)
\end{eqnarray}
Or equivalently in terms of subspaces \begin{equation}\label{eqDP} \cD+\cP=\left((\cD-\cW)+(\cP-\cU)\right)\oplus\left( (\cD+\cW)+(\cP+\cU)\right)\end{equation}

Let $B_+$ be the kernel of the Hilbert subspace $\cB_+$, intersection of the Hilbert subspaces with kernel $(D-W)$ and $(P-U)$, and $B_-$ be the kernel of the Hilbert subspace $\cB_-$,
intersection of the Hilbert subspaces with kernel $(D+W)$ and $(P+U)$. These two kernels are independent, and there sum define a invariant Kre\u{\i}n subspace hence $(B=B_+ -B_-, Q=B_+
+B_-)\in \Gamma_D(\Lambda)$. But $(B,Q)\leq _{\Gamma_D(\Lambda)} (W,D); (B,Q)\leq _{\Gamma_D(\Lambda)} (U,P)$, hence either $$(B,Q)=\alpha (W,D)=\beta (U,P),\quad \alpha >0,\; \beta
>0$$ and $(W,D)$ and $(U,P)$ belong to the same extreme ray or $(B,Q)=0$ and in this case $D$ and $P$ are $\gC$-independent. (Precisely, in equation \ref{eqDP}
$\cD+\cP=(\cD-\cW)\oplus(\cP-\cU)\oplus (\cD+\cW)\oplus(\cP+\cU)$ and the sum $\cD\oplus \cP$ is direct).
\end{proof}

We can now prove theorem \ref{direct}:
\begin{proof}
From the theory of integral of Hilbert subspaces\cite{Thomas05} an integral of Hilbert subspaces is direct if the subspaces are disjoint. But lemma \ref{asM} combined with lemmas
\ref{interM} and \ref{mini} prove that the spaces are disjoint. It follows that the integral of Hilbert subspaces is direct, hence that $\Phi$ is one-to-one. Finally the integral of
Kre\u{\i}n-Hilbert pairs of subspaces is direct.
\end{proof}

\subsection{Extremality and irreducibility}
It is easy to prove that a regularly irreducible Kre\u{\i}n space with kernel $K$ induces extremal pairs of kernels $(K,H)$ for any minimal majorant $H$ of $K$ (by the regular Schur's
lemma \ref{RSLemma}, any self-adjoint projection $P$ is the identity).

However the converse is not true as proves the following example:

\begin{example}
Consider example \ref{ex 1} and $T=J_1$, $\cE=\cK$. Then $(I, J_2)$ is an extremal pair but $\cK$ is reducible ($(I,J_1)$ is not extremal).
\end{example}

%Moreover, it may happen than an integral decomposition of spaces
%does not exists (VOIR exemple 7.3 dans \cite{Ennis79}) NE MARCHE PAS ?
Moreover, it is not clear wether an integral decomposition into irreducible subspaces always exists. The reason is that the set of minimal majorant of $K$ is not bounded as soon as
$\cK$ is not definite, and maximisation procedures on Choquet's conical measures may fail, as in the following (trivial) example:

\begin{example}
Consider the same example with $T=I$. Let $P_n^+=\frac{I-J_n}{2}$ and $P_n^-=\frac{I+J_n}{2}$. then the family $\{J_n,\;n\in \NN^*\}$ is not bounded, so are the families of projections
$\{P_n^+,\;n\in \NN\}$ and $\{P_n^-,\;n\in \NN^*\}$, but they define invariant subspaces and $P_n^+ +P_n^-=I$.
\end{example}

To ensure the regular irreducibility of the pairs of spaces occurring in the decomposition \ref{eq direct integral of spaces} we
make the following hypothesis (FS):\\
Fix an invariant Kre\u{\i}n space $\cK$. Then there exists an isomorphism $J:\cE \longrightarrow \cE$ such that:
\begin{enumerate}
\item $(K, JK)$ is a minimal pair defining
$\cK$
\item and $\forall L$ in $\UU$ (hence verifying equation \ref{UU}),
$JLJ^*=L$.
%\item (Double Commutant symmetry (1)) $$\exists J:\cE \longrightarrow \cE \text{ bijective},\; \forall (w,d) \in \Gamma_D ,\; (w, Jw)\in \Gamma_D$$
%\item (Double Commutant symmetry (2)) $$\exists J:\cE \longrightarrow \cE \text{ bijective},\; \forall (w,d) \in \Gamma_D \text{ and
%minimal},\; (w, Jw)\in \Gamma_D \text{ and minimal}$$
%\item (unicity of the decomposition)\\
%The space $\Gamma_D$ is a lattice.
\end{enumerate}

This hypothesis is strong since we have to know the existence of a special symmetry first. However, if the group $G$ is large enough, it may have a representative $J$ such that
$(K,JK)$ is a minimal pair (this is what we called a fundamental representation), and in this case the second condition is always fulfilled. It is also the case if the algebra
generated by $G$ carries a fundamental symmetry.

\begin{lemma}\label{Jextremal}
Under the hypothesis (FS), any extremal pair $(L,JL)\in \Gamma_D$ defines a regularly irreducible Kre\u{\i}n subspace.
\end{lemma}

\begin{proof}
Suppose the Kre\u{\i}n space $\cL$ with kernel $L$ and Hilbert majorant $\cG$ with kernel $G=JL$ is not regularly irreducible. Then there exists a projection $P$ on $\cL$ such that
$P(\cL)$ is an invariant Kre\u{\i}n space. Its kernel is obviously $PL$, direct calculations give that $PJL$ is positive and a minimal majorant of $PL$ since by hypothesis,
$PJL=PL(J^*)^{-1}=JPLJ^*(J^*)^{-1}=JPL$. We can do the same for the projection $(I-P)$ and it follows that
\begin{equation} (L,JL)=(PL,JPL)+((I-P)L,(J(I-P)L)
\end{equation}
with the three terms in $\Gamma_D$, and $(L,JL)$ is not extremal.
\end{proof}

Now take for minimal pair of kernels $(K,H=JK)$.

\begin{lemma}\label{Jsum}
Suppose $(K,JK)=(W,D)\oplus (U,E)$. Then $D=JW, \; E=JU$.
\end{lemma}

\begin{proof}
Let $P$ be the orthogonal projection in the Hilbert subspace with kernel $JK$ on the subspace with kernel $D$. Then \cite{Schwartz64} $D=PJK$. But direct calculations give also that
$W=PK$ (P is also self-adjoint for the indefinite inner product
induced by $K$). \\
We can now use the hypothesis: \begin{eqnarray*} D&=&PJK =PJKJ^* (J^*)^{-1}\\
&=& PK (J^*)^{-1}=W(J^*)^{-1}\\
&=&JWJ^*(J^*)^{-1}= JW
\end{eqnarray*}
\end{proof}

\begin{lemma}\label{Jint}
In the decomposition \ref{eq_int_kern} $$(K,H)=\int_T (K_t,H_t) dm(t)$$ we have
\begin{equation}
H_t=JK_t\; m-a.s.
 \end{equation}
\end{lemma}

%\begin{proof} TRES MAL DIT
%Let $A$ be the set where $f_t\neq Je_t$. By separability of
%$\cL{(\cF,\cE)}$, this set is measurable and by integrating on $A$
%we get
%$$(w,d)=\int_A (e_t,f_t)
%dm(t)=(v,Jv)-(\int_{T\backslash A} (e_t,f_t) dm(t)$$ But
%
%$$(w,Jw)=\int_A (e_t,Je_t)
%dm(t)= (v,c)=(v,Jv)-\int_{T\backslash A} (e_t,Je_t) dm(t)$$ and also
%
%$$(\int_{T\backslash A} (e_t,f_t)
%dm(t)=\int_{T\backslash A} (e_t,Je_t) dm(t)$$
%
%It follows that $d=Jw$ and by unicity of the measure, $f_t=Je_t$.
%\end{proof}

\begin{proof}
We know that any measure $m$ verifying equation \ref{eq_int_kern} is concentrated on a compact and metrizable (hence separable) set $B$. Let $\{\varphi_n, n\in \NN\}$ be a dense family
of $B$ and define
\begin{equation}N=\{t\in T,\; (K_t,H_t)\neq (K_t,JK_t) \}\end{equation}
\begin{equation}N^+_n=\{t\in N,\; \varphi_n(JK_t-H_t)>0\}\end{equation}
\begin{equation}N^-_n=\{t\in N,\; \varphi_n(JK_t-H_t)<0\}\end{equation}
%\begin{equation}N_n=\{(t,t')\in T^2,\; (e_t+e_{t'}, f_t+f_{t'})=\gamma_{t,t'}+(0,\h_{t,t'}),\; \gamma_t \text{ minimal and }\varphi_n(h_{t,t'})>0\}\end{equation}
Then $N=\bigcup_{n\in \NN} (N^+_n\cup N^-_n)$ by the Hahn-Banach
theorem.\\
By theorem \ref{direct}, we can change the measure such that the following integral is direct:
\begin{equation}
(K,H)=\int_T^{\oplus} (K_t,H_t)dm(t)= \int_{N_n}^{\oplus}(K_t,H_t)dm(t)\oplus \int_{T\backslash N_n}^{\oplus}(K_t,H_t)dm(t)
\end{equation} and by lemma \ref{Jsum} $\int_{N^+_n} (K_t,H_t)
dm(t)$ is of the form $(W,JW)$. But $W=\int_{N^+_n} K_t dm(t)$ hence $JW=\int_{N^+_n} JK_t dm(t)$ and we get $\int_{N^+_n} (JK_t-H_t)dm(t)=0$. It follows that $\int_{N^+_n}
\varphi_n(JK_t-H_t) dm(t)=0$ which implies that $m(N^+_n)=0$ since $\varphi_n(JK_t-H_t)>0$. The same arguments work for $N_n^-$ and finally $m(N)=0$.
\end{proof}

Combining theorem \ref{direct}, lemma \ref{Jint} and lemma \ref{Jextremal} we get:

\begin{theorem}\label{th main}
Let $\cK$ be a invariant Kre\u{\i}n subspace, and suppose that the hypothesis (FS) is verified. Suppose
moreover that the space %$\cL{(\cF,\cE)}$ is conuclear. Then
$\cE$ is weakly complete and conuclear. Then $(\cK,\cH=J\cK)$ admits a direct integral decomposition in terms of irreducible invariant Kre\u{\i}n-Hilbert pairs $(\cK_t,\cH_t=J\cK_t)$:
\begin{equation}
(\cK,\cH)=\int^{\oplus}_T \left(\cK_t ,\cH_t\right) dm(t)
\end{equation}
Moreover, we have an analogue of Parseval's formula:
\begin{equation}
\left[k=\int^{\oplus}_T k_t dm(t),k'=\int^{\oplus}_T k^{'}_t dm(t)\right]_{\cK}= \int_{T} \left[k_t,k^{'}_t\right]_{\cK_t}dm(t)
\end{equation}
\end{theorem}

\begin{example}
Consider example \ref{ex HP}: $X=\RR^3$ and $G=SO(1,2)$ is the associated Lorentz group, $J=\left(
                                                 \begin{array}{ccc}
                                                   1 & 0 & 0 \\
                                                   0 & -1 & 0 \\
                                                   0 & 0 & -1 \\
                                                 \end{array}
                                               \right)$
is the Minkowski metric operator.\\

Note that $\cE=\CC^{X}$ is a nuclear and Frechet space (as a product of nuclear and Frechet spaces \cite{Treves67}), hence conuclear
\cite{Groth55}, \cite{Schwartz73}.\\
%Note that $\cE=\CC^{X}$ is complete and conuclear. Indeed, a product
%of nuclear spaces is nuclear (proposition 50.1 in \cite{Treves67})
%and therefore $\cE$ is nuclear. But it is barreled as the product of
%barreled spaces, and reflexive. It follows that its dual is also
%barreled and nuclear, hence $\cE$ is conuclear. In fact $\CC^{X}$ is
%a nuclear Frechet space, hence conuclear \cite{Groth55},
%\cite{Schwartz73}.\\

The pair $K(x,y)=\exp{\langle Jx,y\rangle_{\RR^3}}$, $H(x,y)=\exp{\langle x,y\rangle_{\RR^3}}$ is a fundamental
%Krein-Hilbert
pair of kernels, $K(x,y)$ is invariant under $G$ and for all $g\in G$ $\tau_g(\cH)\subset \cH$. Moreover $K=\tau_J H$ with $J$ in $G$, and the representation  is fundamental. Finally
the hypothesis of theorem \ref{th main} are fulfilled and $(\cK,\cH=\tau_ J\cK)$ admits a direct integral decomposition in terms of irreducible invariant Kre\u{\i}n-Hilbert pairs
$(\cK_t,\cH_t=\tau_J\cK_t)$.

The decomposition is as follows:
\begin{equation}
K_n(x,y)=\frac{(\langle Jx,y\rangle_{\RR^3})^n}{n!},\qquad H_n(x,y)=JK_n(x,y)=\frac{(\langle x,y\rangle_{\RR^3})^n}{n!}\end{equation} and
\begin{equation}
(\cK,\cH)=\bigoplus (\cK_n,\cH_n)
\end{equation}

Note that each Kre\u{\i}n space $\cK_n$ of homogeneous polynomials of degree $n$ with kernel $K_n(x,y)$ is regularly irreducible by theorem \ref{th_class1}.
\end{example}

%\section*{Conclusion and comments}

\bibliographystyle{amsplain}
%\bibliography{Bibliototale}

%\makebibliography

\end{document}